\documentclass[12pt,reqno]{amsart}
\usepackage{times}

\newtheorem{thm}{Theorem}
\newtheorem{lemma}[thm]{Lemma}
\newtheorem{cor}[thm]{Corollary}

\newcommand{\cn}{{\mathbb C}^n}
\newcommand{\D}{{\mathbb D}}
\newcommand{\ind}{\int_{\D}}
\newcommand{\bn}{{\mathbb B}_n}
\newcommand{\inb}{\int_{\bn}}
\newcommand{\sn}{{\mathbb S}_n}
\newcommand{\ins}{\int_{\sn}}

\newcommand{\radd}{(1-|z|^2)|Rf(z)|}
\newcommand{\grad}{(1-|z|^2)|\nabla f(z)|}
\newcommand{\invd}{|\widetilde\nabla f(z)|}
\newcommand{\aut}{{\rm Aut}(\bn)}

\begin{document}

\title[Bergman Spaces]
{New Characterizations of Bergman Spaces}

\author{Miroslav Pavlovi\'c and Kehe Zhu}
\address{Miroslav Pavlovi\'c\\
              Matemati\v cki Fakultet\\
              Studentski Trg 16\\
              11001 Belgrade, P.P. 550\\
              Serbia}
\email{pavlovic@matf.bg.ac.yu}
\address{Kehe Zhu\\
              Department of Mathematics\\
              SUNY\\
              Albany, NY 12222, USA}
\email{kzhu@math.albany.edu}
\date{January 30, 2006}
\thanks{The first author is supported in part by MNZZS Grant ON144010 and
the second author is partially supported by the National Science Foundation}

\begin{abstract}
We obtain several new characterizations for the standard weighted Bergman spaces 
$A^p_\alpha$ on the unit ball of $\cn$ in terms of the radial derivative, the holomorphic 
gradient, and the invariant gradient.
\end{abstract}

\maketitle

\section{Introduction}

Let $\bn$ be the open unit ball in $\cn$. For $\alpha>-1$ let $$dv_\alpha(z)=c_\alpha(1-|z|^2)^\alpha\,dv(z),$$
where $dv$ is the normalized volume measure on $\bn$ and $c_\alpha$
is a positive constant making $dv_\alpha$ a probability measure. For $0<p<\infty$ 
the weighted Bergman space $A^p_\alpha$ consists of holomorphic functions in
$L^p(\bn,dv_\alpha)$. Thus
$$A^p_\alpha= H(\bn)\cap L^p(\bn,dv_\alpha),$$
where $H(\bn)$ is the space of all holomorphic functions in $\bn$.

For $f\in H(\bn)$ and $z=(z_1,\cdots,z_n)\in\bn$ we define
$$Rf(z)=\sum_{k=1}^nz_k\frac{\partial f}{\partial z_k}(z)$$
and call it the radial derivative of $f$ at $z$. The complex gradient of $f$ at $z$
is defined as
$$|\nabla f(z)|=\left[\sum_{k=1}^n\left|\frac{\partial f}{\partial z_k}(z)\right|^2
\right]^{1/2}.$$

Let $\aut$ denote the automorphism group of $\bn$. Thus $\aut$ consists of all
bijective holomorphic functions $\varphi:\bn\to\bn$. It is well known that $\aut$ is
generated by two types of maps: unitaries and symmetries. The unitaries are simiply
the $n\times n$ unitary matrices considered as mappings from $\bn$ to $\bn$. For
any point $a\in\bn$ there exists a unique map $\varphi_a\in\aut$ with the following
properties: $\varphi_a(0)=a$, $\varphi_a(a)=0$, and $\varphi_a\circ\varphi_a(z)=z$
for all $z\in\D$. Such a mapping $\varphi_a$ is called a symmetry. Because of the
property $\varphi_a\circ\varphi_a(z)=z$ it is also natural to call $\varphi_a$ an
involution or an involutive automorphism. See \cite{rudin} and \cite{zhu} for more 
information about the automorphism group of $\bn$.

If $f\in H(\bn)$, we define
$$|\widetilde\nabla f(z)|=|\nabla(f\circ\varphi_z)(0)|,\qquad z\in\bn.$$
It can be checked that 
$$|\widetilde\nabla(f\circ\varphi)|=|(\widetilde\nabla f)\circ
\varphi|,\qquad \varphi\in\aut.$$
So $|\widetilde\nabla f(z)|$ is called the invariant gradient of $f$ at $z$. 
See \cite{zhu} for more information about the invariant gradient.

When $n=1$, the unit ball ${\mathbb B}_1$ is usually called the unit disk and we
denote it by $\D$ instead. In this case, we clearly have
$$Rf(z)=zf(z),\quad |\nabla f(z)|=|f'(z)|,\quad |\widetilde\nabla f(z)|=(1-|z|^2)|f'(z)|.$$
In particular, the functions
\begin{equation}
(1-|z|^2)|Rf(z)|,\quad (1-|z|^2)|\nabla f(z)|,\quad |\widetilde\nabla f(z)|,
\label{eq1}
\end{equation}
have exactly the same boundary behavior on the unit disk $\D$. In higher dimensions,
the three functions above no longer have the same boundary behavior; see Section 2.3
and Chapter 7 in \cite{zhu}. However, when integrated against the weighted volume
measures $dv_\alpha$, not only do these differential-based functions exhibit the same
behavior, they also behave the same as the original function $f(z)$, as the following
result (see Theorem 2.16 of \cite{zhu}) demonstrates.

\begin{thm}
Suppose $p>0$, $\alpha>-1$, and $f\in H(\bn)$. Then the following conditions are
equivalent.
\begin{enumerate}
\item[(a)] $f\in A^p_\alpha$, that is, $f\in L^p(\bn,dv_\alpha)$.
\item[(b)] The function $f_1(z)=(1-|z|^2)|Rf(z)|$ belongs to $L^p(\bn,dv_\alpha)$.
\item[(c)] The function $f_2(z)=(1-|z|^2)|\nabla f(z)|$ belongs to $L^p(\bn,dv_\alpha)$.
\item[(d)] The function $f_3(z)=|\widetilde\nabla f(z)|$ belongs to $L^p(\bn,dv_\alpha)$.
\end{enumerate}
Moreover, the quantities
$$|f(0)|^p+\inb|f_1|^p\,dv_\alpha,\ 
|f(0)|^p+\inb|f_2|^p\,dv_\alpha,\  |f(0)|^p+\inb|f_3|^p\,dv_\alpha,$$
are all comparable to
$$\inb|f(z)|^p\,dv_\alpha(z)$$
whenever $f$ is holomorphic in $\bn$.
\label{1}
\end{thm}

The purpose of this paper is to explore the above ideas further. We show
that the integral behavior of the functions
$$|f(z)|,\quad (1-|z|^2)|Rf(z)|,\quad (1-|z|^2)|\nabla f(z)|,\quad
|\widetilde\nabla f(z)|,$$ 
is the same in a much stronger sense. More specifically, when integrating over
the unit ball with respect to weighted volume measures, we can write
$|f(z)|^p=|f(z)|^{p-q}|f(z)|^q$ and can replace $|f(z)|$ in the second factor
by any one of the functions in (\ref{eq1}). We state our main result as follows.

\begin{thm}
Suppose $p>0$, $\alpha>-1$, $0<q<p+2$, and $f\in H(\bn)$. Then the following 
conditions are equivalent.
\begin{enumerate}
\item[(a)] $f\in A^p_\alpha$, that is, $I_1(f)<\infty$, where
$$I_1(f)=\inb|f(z)|^p\,dv_\alpha(z).$$
\item[(b)] $I_2(f)<\infty$, where
$$I_2(f)=\inb|f(z)|^{p-q}\left[\radd\right]^q\,dv_\alpha(z).$$
\item[(c)] $I_3(f)<\infty$, where
$$I_3(f)=\inb|f(z)|^{p-q}\left[\grad\right]^q\,dv_\alpha(z).$$
\item[(d)] $I_4(f)<\infty$, where
$$I_4(f)=\inb|f(z)|^{p-q}\invd^q\,dv_\alpha(z).$$
\end{enumerate}
Furthermore, the quantities
$$I_1(f),\quad |f(0)|^p+I_2(f),\quad |f(0)|^p+I_3(f),\quad |f(0)|^p+I_4(f),$$
are comparable for $f\in H(\bn)$.
\label{2}
\end{thm}

We will show by a simple example that the range $0<q<p+2$ is best possible.

Throughout the paper we use $C$ to denote a positive constant,
indepedent of $f$ and $z$, whose value may vary from one occurence to another.

\section{The case $0<q\le p$}

The proof of Theorem~\ref{2} requires different methods for the two
cases $0<q\le p$ and $p<q<p+2$. This section deals with the case $0<q\le p$;
the other case is considered in the next section.

The case $q=p$ is of course just Theorem~\ref{1}. Our proof of Theorem~\ref{2} in 
the case $0<q<p$ is based on several technical lemmas that are known to experts. 
We include them here for the non-expert and for convenience of reference. We begin 
with the following embedding theorem for Bergman spaces.

\begin{lemma}
Suppose $0<p\le1$, $\alpha>-1$, and
$$\beta=\frac{n+1+\alpha}p-(n+1).$$
There exists a constant $C>0$ such that
$$\inb|f(z)|\,dv_\beta(z)\le C\left[\inb|f(z)|^p\,dv_\alpha(z)\right]^{1/p}$$
for all $f\in H(\bn)$.
\label{3}
\end{lemma}

\begin{proof}
See Lemma 2.15 of \cite{zhu}.
\end{proof}

We will also need the following boundedness criterion for a class of
integral operators on $\bn$.

\begin{lemma}
For real $a$ and $b$ consider the integral operator $T=T_{a,b}$ defined by
$$Tf(z)=(1-|z|^2)^a\inb\frac{(1-|w|^2)^b}{|1-\langle z,w\rangle|^{n+1+a+b}}
f(w)\,dv(w),$$
where
$$\langle z,w\rangle=\sum_{k=1}^nz_k\overline w_k$$
for $z=(z_1,\cdots,z_n)$ and $w=(w_1,\cdots,w_n)$ in $\bn$. If $p\ge1$, then
$T$ is bounded on $L^p(\bn,dv_\alpha)$ if and only if the inequalities
$$-pa<\alpha+1<p(b+1)$$
hold.
\label{4}
\end{lemma}

\begin{proof}
See Theorem 2.10 of \cite{zhu}.
\end{proof}

The following result compares the various derivatives that we use for
a holomorphic function in $\bn$.

\begin{lemma}
If $f\in H(\bn)$, then
$$\invd^2=(1-|z|^2)(|\nabla f(z)|^2-|Rf(z)|^2).$$
Moreover,
$$(1-|z|^2)|Rf(z)|\le(1-|z|^2)|\nabla f(z)|\le|\widetilde\nabla f(z)|$$
for all $z\in\bn$.
\label{5}
\end{lemma}

\begin{proof}
See Lemmas 2.13 and 2.14 of \cite{zhu}.
\end{proof}

We will need the following well-known reproducing formula for holomorphic
functions in $\bn$.

\begin{lemma}
If $\alpha>-1$ and $f\in A^1_\alpha$, then
$$f(z)=\inb\frac{f(w)\,dv_\alpha(w)}{(1-\langle z,w\rangle)^{n+1+\alpha}}$$
for all $z\in\bn$.
\label{6}
\end{lemma}

\begin{proof}
See Theorem 2.2 of \cite{zhu}.
\end{proof}

The following integral estimate is standard in the theory of Bergman spaces
and has proved to be very useful in many different situations.

\begin{lemma}
Suppose $\alpha>-1$  and $t>0$. Then there exists a constant $C>0$ such that
$$\inb\frac{dv_\alpha(w)}{|1-\langle z,w\rangle|^{n+1+\alpha+t}}\le
\frac C{(1-|z|^2)^t}$$
for all $z\in\bn$.
\label{7}
\end{lemma}

\begin{proof}
See Proposition 1.4.10 of \cite{rudin} or Theorem 1.12 of \cite{zhu}.
\end{proof}

We now begin the proof of Theorem~\ref{2} under the assumption that $0<q<p$. 
In this case, the numbers $r=p/(p-q)$ and $s=p/q$ satisfy $r>1$, $s>1$, and
$1/r+1/s=1$. So we can apply H\"older's inequality to the integral $I_4(f)$ to obtain
\begin{equation}
I_4(f)\le\left[\inb|f(z)|^p\,dA_\alpha(z)\right]^{\frac1r}
\left[\inb|\widetilde\nabla f(z)|^p\,dv_\alpha(z)\right]^{\frac1s}.
\label{eq2}
\end{equation}
By Theorem~\ref{1}, there exists a positive constant $C>0$, independent of $f$,
such that
$$\inb|\widetilde\nabla f(z)|^p\,dv_\alpha(z)\le C\inb|f(z)|^p\,dv_\alpha(z).$$
Combining this with (\ref{eq2}), we see that the integral  $I_4(f)$ is dominated by
$I_1(f)$.

According to Lemma~\ref{5}, we have $I_2(f)\le I_3(f)\le I_4(f)$. So it remains for
us to show that $I_1(f)$ is finite whenever $I_2(f)$ is finite. We do this in two steps.

First, we assume that $p=qN$ for some integer $N>1$. In this case, the function
$f(z)^{p/q}$ is well-defined and holomorphic in $\bn$. Moreover,
$$R\left[f(z)^{\frac pq}\right]=\frac pq\,f(z)^{\frac pq-1}Rf(z).$$
Let $\beta$ be a sufficiently large (to be specified later) positive integer and
apply Lemma~\ref{6} to write
$$R\left[f(z)^{\frac pq}\right]=\frac pq\inb\frac{f(w)^{\frac pq-1}Rf(w)\,dv_\beta(w)}
{(1-\langle z,w\rangle)^{n+1+\beta}},\qquad z\in\bn.$$
Since the function $f(w)^{(p/q)-1}Rf(w)$ vanishes at the origin, we can also write
$$R\left[f(z)^{\frac pq}\right]=\frac pq\inb\left[\frac1{(1-\langle z,w\rangle)^{n+1
+\beta}}-1\right]f(w)^{\frac pq-1}Rf(w)\,dv_\beta(w).$$
Integrating the above equation, we obtain
$$f(z)^{\frac pq}-f(0)^{\frac pq}=\int_0^1\!Rf^{\frac pq}(tz)\,\frac{dt}t
=\inb\!H(z,w)f(w)^{\frac pq-1}Rf(w)\,dv_\beta(w),$$
where
$$H(z,w)=\frac pq\int_0^1\frac{1-(1-t\langle z,w\rangle)^{n+1+\beta}}
{(1-t\langle z,w\rangle)^{n+1+\beta}}\,\frac{dt}t.$$
Expand the numerator in the integrand above by the binomial formula and then
evaluate the integral term by term. We obtain a positive constant $C>0$ such that
$$|H(z,w)|\le\frac C{|1-\langle z,w\rangle|^{n+\beta}}$$
for all $z$ and $w$ in $\bn$. It follows that
\begin{equation}
\left|f(z)^{\frac pq}-f(0)^{\frac pq}\right|\le C\inb\frac{|f(w)|^{\frac pq-1}|Rf(w)|
\,dv_\beta(w)}{|1-\langle z,w\rangle|^{n+\beta}}
\label{eq3}
\end{equation}
for all $z\in\bn$.

If $q\ge1$, then we rewrite (\ref{eq3}) as
\begin{equation}
\left|f(z)^{\frac pq}-f(0)^{\frac pq}\right|\le C\inb g(w)
\frac{(1-|w|^2)^{\beta-1}\,dv(w)}{|1-\langle z,w\rangle|^{n+1+\beta-1}},
\label{eq4}
\end{equation}
where
$$g(w)=|f(w)|^{\frac pq-1}(1-|w|^2)|Rf(w)|.$$
By Lemma~\ref{4}, the integral operator
$$Tg(z)=\inb g(w)\,\frac{(1-|w|^2)^{\beta-1}\,dv(w)}{|1-\langle z,
w\rangle|^{n+1+\beta-1}}$$
is bounded on $L^q(\bn,dv_\alpha)$, because we can choose the positive integer 
$\beta$ to satisfy $\alpha+1<q\beta$. Combining this with (\ref{eq4}), we obtain
a positive constant $C$, independent of $f$, such that
$$\inb\left|f^{\frac pq}-f(0)^{\frac pq}\right|^qdv_\alpha\le 
C\inb|f(z)|^{p-q}\left[(1-|z|^2)|Rf(z)|\right]^qdv_\alpha(z).$$
This clearly shows that there exists a positive constant $C>0$, independent of $f$,
such that
$$I_1(f)\le C\left[|f(0)|^p+I_2(f)\right]$$
for all $f\in H(\bn)$.

If $0<q<1$, we rewrite (\ref{eq3}) as
\begin{equation}
\left|f(z)^{\frac pq}-f(0)^{\frac pq}\right|\le C\inb\left|\frac{f(w)^{\frac pq-1}Rf(w)}
{(1-\langle w,z\rangle)^{n+\beta}}\right|(1-|w|^2)^\beta\,dv(w).
\label{eq5}
\end{equation}
We also write
$$\beta=\frac{n+1+\gamma}q-(n+1),$$
and choose $\beta$ to be large enough so that $\gamma>-1$. We then apply 
Lemma~\ref{3} to the right-hand side of (\ref{eq5}) to obtain
$$\left|f(z)^{\frac pq}-f(0)^{\frac pq}\right|\le C\left[\inb\left|
\frac{f(w)^{\frac pq-1}Rf(w)}{(1-\langle z,w\rangle)^{n+\beta}}\right|^q
\,dv_\gamma(w)\right]^{\frac 1q},$$
where $C$ is a positive constant independent of $f$. Take the $q$th power on
both sides, integrate over $\bn$ with respect to $dv_\alpha$, and apply Fubini's
theorem. We see that the integral
$$\inb\left|f(z)^{\frac pq}-f(0)^{\frac pq}\right|^q\,dv_\alpha$$
is dominated by the integral
$$\inb|f(w)|^{p-q}|Rf(w)|^q\,dv_\gamma(w)
\inb\frac{dv_\alpha(z)}{|1-\langle z,w\rangle|^{q(n+\beta)}}.$$
If $\beta$ is large enough so that
$$q(n+\beta)>n+1+\alpha,$$
then by Lemma~\ref{7}, there exists a positive constant $C$ such that
$$\inb\frac{dv_\alpha(z)}{|1-\langle z,w\rangle|^{q(n+\beta)}}\le\frac C{(1-
|w|^2)^{q(n+\beta)-(n+1+\alpha)}}$$
for all $w\in\bn$. An easy calculation shows that
$$q(n+\beta)-(n+1+\alpha)=\gamma-(q+\alpha).$$
It follows that
$$\inb\left|f^{\frac pq}-f(0)^{\frac pq}\right|^qdv_\alpha\le C\inb|f(z)|^{p-q}\left[(1-|z|^2)|Rf(z)|\right]^qdv_\alpha(z),$$
where $C$ is a positive constant independent of $f$. This easily implies that
$$I_1(f)\le C\left[|f(0)|^p+I_2(f)\right]$$
for another positive constant $C$ that is independent of $f$.

Thus we have proved that the integral $I_1(f)$ is dominated by $|f(0)|^p+I_2(f)$
under the additional assumption that $p=qN$, where $N>1$ is a positive integer.

In the general case $0<q<p$, we choose a positive integer $N$ such that $Nq>p$
and define two positive numbers $r$ and $s$ by
$$r=\frac{Nq}p,\qquad \frac1r+\frac1s=1.$$
By the special case that we have already proved, there exists a constant $C>0$,
independent of $f$, such that
$$I_1(f)\le C\left[|f(0)|^p+\inb\left[|f(z)|^{-1}(1-|z|^2)|Rf(z)|\right]^{p/N}
|f(z)|^p\,dv_\alpha(z)\right].$$
By an approximation argument we may assume that $I_1(f)$ is finite (note that 
we are trying to prove the stronger conclusion that $I_1(f)$ is dominated by 
$|f(0)|^p+I_2(f)$). By H\"older's inequality, the integral on the right-hand side 
above does not exceed
$$\left[\inb\!\left[|f(z)|^{-1}(1-|z|^2)|Rf(z)|\right]^{rp/N}|f(z)|^p\,dv_\alpha(z)
\right]^{\frac1r}\!\left[\inb\!|f|^p\,dv_\alpha\right]^{\frac1s}.$$
It follows that
$$I_1(f)\le C\left[|f(0)|^p+I_2(f)^{\frac1r}I_1(f)^{\frac 1s}\right].$$
From this we easily deduce that $I_1(f)$ is dominated by $|f(0|^p+I_2(f)$.
In fact, this is obvious if $f(0)=0$. Otherwise, we may use homogeneity to assume 
that $f(0)=1$. In this case, we also have $I_1(f)\ge1$, so dividing both sides of the 
above inequality by $I_1(f)^{1/s}$ yields
$$I_1(f)^{\frac1r}\le C\left[\frac1{I_1(f)^{1/s}}+I_2(f)^{\frac1r}\right]
\le C\left[1+I_2(f)^{\frac1r}\right].$$
This clearly implies that
$$I_1(f)\le C\left[1+I_2(f)\right]=C\left[|f(0)|^p+I_2(f)\right]$$
for some other positive constant independent of $f$.
This completes the proof of Theorem~\ref{2} in the case $0<q\le p$.

\section{The case $p<q<p+2$}

This section is devoted to the proof of Theorem~\ref{2} in the case $p<q<p+2$.

It follows from Theorem~\ref{1} that there exists a small positive constant $c$
such that
\begin{eqnarray*}
cI_1(f)-|f(0)|^p&\le& \inb(1-|z|^2)^p|Rf(z)|^p\,dv_\alpha(z)\\
&=& \inb(1-|z|^2)^p|Rf(z)|^p|f(z)|^a|f(z)|^{-a}\,dv_\alpha(z),
\end{eqnarray*}
where $a=p(p-q)/q$. Let
$$r=\frac qp,\qquad s=\frac q{q-p}.$$
When $p<q$, we have $r>1$, $s>1$, and $1/r+1/s=1$. An application of
H\"older's inequality shows that $cI_1(f)-|f(0)|^p$ does not exceed
$$\left[\inb(1-|z|^2)^q|Rf(z)|^q|f(z)|^{p-q}\,dv_\alpha(z)
\right]^{\frac1r}\left[\inb|f(z)|^p\,dv_\alpha(z)\right]^{\frac 1s}.$$
Therefore,
$$cI_1(f)\le|f(0)|^p+I_2(f)^{\frac1r}I_1(f)^{\frac1s}.$$
From this we easily deduce that 
$$I_1(f)\le C\left[|f(0)|^p+I_2(f)\right]$$
for some positive constant $C$ independent of $f$; see the last paragraph of
the previous section.

Once again, Lemma~\ref{5} tells us that $I_2(f)\le I_3(f)\le I_4(f)$. So it
remains for us to show that the integral $I_4(f)$ is dominated by $I_1(f)$.
This will require several technical lemmas again.

We begin with the following well-known estimate for the Bergman kernel on
pseudo-hyperbolic balls.

\begin{lemma}
Suppose $\rho\in(0,1)$. Then there exists a positive constant $C$ (independent
of $z$ and $w$) such that
$$C^{-1}(1-|z|^2)\le|1-\langle z,w\rangle|\le C(1-|w|^2)$$
for all $z$ and $w$ in $\bn$ satisfying $|\varphi_z(w)|<\rho$. Moreover, if 
$$D(z,\rho)=\{w\in\bn:|\varphi_z(w)|<\rho\}$$
is a pseudo-hyperbolic ball, then its Euclidean volume satisfies
$$C^{-1}(1-|z|^2)^{n+1}\le v(D(z,\rho))\le C(1-|z|^2)^{n+1}.$$
\label{8}
\end{lemma}

\begin{proof}
See Lemmas 1.23 and 2.20 of \cite{zhu}.
\end{proof}

Note that, by symmetry, the positions of $z$ and $w$ can be interchanged in the
first set of inequalities of Lemma~\ref{8}.

The key to the remaining proof of Theorem~\ref{2} is the following well-known
special case of $q=2$.

\begin{lemma}
For every $p>0$ there exists a positive constant $C$ such that
$$\inb|f(z)|^p\,dv(z)\le C\left[|f(0)|^p+\inb|f(z)|^{p-2}\invd^2\,dv(z)\right]$$
and
$$|f(0)|^p+\inb|f(z)|^{p-2}\invd^2\,dv(z)\le C\inb|f(z)|^p\,dv(z)$$
for all $f\in H(\bn)$.
\label{9}
\end{lemma}

\begin{proof}
See \cite{oyz}.
\end{proof}

In the general case, we first prove the following weaker version.

\begin{lemma}
Suppose $p>0$, $0<q<p+2$, and $\alpha>-1$. There exists
a positive constant $C$ (independent of $f$) such that
$$\int_{|z|<1/4}|f(z)|^{p-q}\invd^q\,dv_\alpha(z)\le C\int_{|z|<3/4}|f(z)|^p
\,dv_\alpha(z)$$
for all $f\in H(\bn)$.
\label{10}
\end{lemma}

\begin{proof}
If $0<q\le p$, the desired estimate follows from the well-known fact that 
point-evaluations (of any form of the derivative) on a compact subset of $|z|<3/4$ 
are uniformly bounded linear functionals on the Bergman spaces of the ball 
$|z|<3/4$; see Lemma 2.4 of \cite{zhu} for example.

So we assume that $p<q<p+2$. In this case, we have $1<2/(q-p)$. Fix
$r\in(1,2/(q-p))$, sufficiently close to $2/(q-p)$, so that $q-\lambda>0$,
where $\lambda=2/r\in(q-p,2)$.

If $f$ is a unit vector in $H^\infty(\bn)$, then there exists a constant $C>0$,
independent of $f$, such that $|\nabla f(0)|\le C$. Replacing $f$ by $f\circ\varphi_z$,
we obtain $\invd\le C$ for all $z\in\bn$. It follows from this and H\"older's
inequality that the integral
$$I(f)=\int_{|z|<1/2}|f(z)|^{p-q}\invd^q\,dv(z)$$
satisfies
\begin{eqnarray*}
I(f)&=&\int_{|z|<1/2}|f(z)|^{p-q}\invd^\lambda\invd^{q-\lambda}\,dv(z)\\
&\le& C^{q-\lambda}\int_{|z|<1/2}|f(z)|^{p-q}\invd^\lambda\,dv(z)\\
&\le& C^{q-\lambda}\left[\int_{|z|<1/2}|f(z)|^{r(p-q)}\invd^{r\lambda}
\,dv(z)\right]^{\frac1r}\\
&\le& C^{q-\lambda}\left[\inb|f(z)|^{r(p-q)}\invd^{r\lambda}
\,dv(z)\right]^{\frac1r}\\
&=& C^{q-\lambda}\left[\inb|f(z)|^{r(p-q)+2-2}\invd^2\,dv(z)\right]^{\frac1r}.
\end{eqnarray*}
By Lemma~\ref{9}, there exists a positive constant $C$, independent of $f$,
such that
$$I(f)\le C\left[\inb|f(z)|^{r(p-q)+2}\,dv(z)\right]^{\frac1r}\le C$$
for all unit vectors $f$ of $H^\infty(\bn)$. Here we used the assumption that
$r(p-q)+2>0$, which is equivalent to $r<2/(q-p)$. If $f$ is an arbitrary function
in $H^\infty(\bn)$, then replacing $f$ by $f/\|f\|_\infty$ in $I(f)\le C$
leads to
\begin{equation}
\int_{|z|<1/2}|f(z)|^{p-q}\invd^q\,dv(z)\le C\|f\|_\infty^p,
\label{eq6}
\end{equation}
where
$$\|f\|_\infty=\sup\{|f(z)|:z\in\bn\}.$$

It is easy to see that $\invd$ and $|\nabla f(z)|$ are comparable on any
compact subset of $\bn$. In fact, it follows from Lemma~\ref{5} that
$$(1-|z|^2)|\nabla f(z)|\le\invd\le|\nabla f(z)|,$$
which shows that $\invd$ and $|\nabla f(z)$ are comparable on any compact
subset of $\bn$.

Now suppose $f$ is any holomorphic function in $\bn$. We replace
$f(z)$ in (\ref{eq6}) by $f(z/2)$, use the conclusion of the previous paragraph,
and make the change of variables $w=z/2$. Then there exists a positive constant
$C$, independent of $f$, such that
$$\int_{|z|<1/4}|f(z)|^{p-q}\invd^q\,dv(z)\le C\sup\{|f(z)|^p:|z|\le1/2\}.$$
Since point-evaluations in $|z|\le1/2$ are uniformly bounded on Bergman
spaces of the ball $|z|<3/4$, there exists a positive constant $C$, independent of $f$, 
such that
$$\int_{|z|<1/4}|f(z)|^{p-q}\invd^q\,dv(z)\le C\int_{|z|<3/4}|f(z)|^p\,dv(z).$$
Since $(1-|z|^2)^\alpha$ is comparable to a positive constant whenever $z$ is
restricted to a compact subset of $\bn$, we obtain a positive constant $C$,
independent of $f$, such that
$$\int_{|z|<1/4}|f(z)|^{p-q}\invd^q\,dv_\alpha(z)\le 
C\int_{|z|<3/4}|f(z)|^p\,dv_\alpha(z).$$
This completes the proof of Lemma~\ref{10}.
\end{proof}

We now use Lemma~\ref{10} to show that the integral $I_4(f)$ is dominated by $I_1(f)$. 
This part of the proof works for the full range $0<q<p+2$.

Replace $f$ by $f\circ\varphi_w$ in Lemma~\ref{10}, where $w$ is an arbitrary
point in $\bn$, and use the M\"obius invariance of $\widetilde\nabla f$. Then
the integrals
$$\int_{|z|<1/4}|f(\varphi_w(z))|^{p-q}|(\widetilde\nabla f)(\varphi_w(z))|^q
\,dv_\alpha(z)$$
are uniformly (with respecto to $w$) dominated by the integrals
$$\int_{|z|<3/4}|f(\varphi_w(z))|^p\,dv_\alpha(z).$$
Making the change of variables $z\mapsto\varphi_w(z)$ in the above integrals,
we see that the integrals
$$\int_{|\varphi_w(z)|<1/4}|f(z)|^{p-q}\invd^q\frac{(1-|w|^2)^{n+1+\alpha}}
{|1-\langle z,w\rangle|^{2(n+1+\alpha)}}\,dv_\alpha(z)$$
are uniformly (with respect to $w$) dominated by the integrals
$$\int_{|\varphi_w(z)|<3/4}|f(z)|^p\frac{(1-|w|^2)^{n+1+\alpha}}
{|1-\langle z,w\rangle|^{2(n+1+\alpha)}}\,dv_\alpha(z).$$
According to Lemma~\ref{8}, for $|\varphi_w(z)|<3/4$ (hence for 
$|\varphi_w(z)|<1/4$ as well) we have
$$1-|w|^2\sim1-|z|^2\sim|1-\langle z,w\rangle|.$$
It follows that there exists another positive constant $C$, independent of $f$ and $w$,
such that
$$\int_{|\varphi_w(z)|<1/4}|f(z)|^{p-q}\invd^q\,dv_\alpha(z)
\le C\int_{|\varphi_w(z)|<3/4}|f(z)|^p\,dv_\alpha(z)$$
for all $f\in H(\bn)$. Integrate the above inequality
over $\bn$ with respect to the M\"obius invariant measure
$$d\tau(w)=\frac{dv(w)}{(1-|w|^2)^{n+1}}.$$
We see that the integral
\begin{equation}
\inb d\tau(w)\int_{|\varphi_z(w)|<1/4}|f(z)|^{p-q}|\invd^q\,dv_\alpha(z)
\label{eq7}
\end{equation}
is dominated by the integral
\begin{equation}
\inb d\tau(w)\int_{|\varphi_z(w)|<3/4}|f(z)|^p\,dv_\alpha(z).
\label{eq8}
\end{equation}

By Fubini's theorem, the integral in (\ref{eq7}) equals
$$\inb|f(z)|^{p-q}\invd^q\,dv_\alpha(z)\int_{|\varphi_w(z)|<1/4}\,d\tau(w).$$
Similarly, the integral in (\ref{eq8}) equals
$$\inb|f(z)|^p\,dv_\alpha(z)\int_{|\varphi_w(z)|<3/4}\,d\tau(w).$$
For any fixed radius $\rho\in(0,1)$, it follows from Lemma~\ref{8} that the integral
$$\int_{|\varphi_w(z)|<\rho}\,d\tau(w)$$
is comparable to a positive constant. Combining these conclusions with (\ref{eq7}) 
and (\ref{eq8}), we obtain another positive constant $C$, independent of $f$, such that
$$\inb|f(z)|^{p-q}\invd^q\,dv_\alpha(z)\le C\inb|f(z)|^p\,dv_\alpha(z)$$
for all $f\in H(\bn)$. This shows that the integral $I_4(f)$ is always dominated by
$I_1(f)$. The proof of Theorem~\ref{2} is now complete.

\section{Further Remarks}

An immediate consequence of Theorem~\ref{2} is the following characterization
of Bergman spaces in terms of the familiar first order partial derivatives.

\begin{cor}
Suppose $p>0$, $0<q<p+2$, $\alpha>-1$, and $f$ is holomorphic in $\bn$. Then 
$f\in A^p_\alpha$ if and only if
\begin{equation}
\inb|f(z)|^{p-q}\left[(1-|z|^2)\left|\frac{\partial f}{\partial z_k}(z)\right|\right]^q
\,dv_\alpha(z)<\infty
\label{eq9}
\end{equation}
for all $1\le k\le n$.
\label{11}
\end{cor}

\begin{proof}
It is clear from the definition of $|\nabla f(z)|$ that for a holomorphic
function $f$ in $\bn$, condition (c) in Theorem~\ref{2} is equivalent to
the condition in (\ref{eq9}).
\end{proof}

Finally we use an example to show that the range $0<q<p+2$ in Theorem~\ref{2}
is best possible. Simply take $f(z)=z_1$. Then on the compact set $|z|\le1/2$,
we have $\invd\sim|\nabla f(z)|=1$. It follows that
\begin{eqnarray*}
\int_{|z|<1/2}|f(z)|^{p-q}\invd^q\,dv_\alpha(z)
&\sim&\int_{|z|<1/2}|f(z)|^{p-q}\,dv_\alpha(z)\\
&=&\int_{|z|<1/2}|z_1|^{p-q}\,dv_\alpha(z).
\end{eqnarray*}
By integration in polar coordinates (see Lemma 1.8 of \cite{zhu} for example),
the last integral above is comparable to
$$\int_0^{1/2}r^{2n-1+p-q}\,dr\ins|\zeta_1|^{p-q}\,d\sigma(\zeta).$$
If $q\ge p+2$, the product above is always infinite. In fact, if $n=1$, then
$$\int_0^{1/2}r^{2n-1+p-q}\,dr=\infty;$$
if $n\ge2$, then by a well-known formula for evaluating integrals of functions of fewer 
variables on the unit sphere (see Lemma 1.9 of \cite{zhu} for example), we have
$$\ins|\zeta_1|^{p-q}\,d\sigma(\zeta)=c\ind|w|^{p-q}(1-|w|^2)^{n-2}\,dA(w)
=\infty,$$
where $c$ is a positive constant and $dA$ is area measure on the unit disk $\D$. 
This shows that the range $q<p+2$ is best possible in Theorem~\ref{2} as well
as in Lemma~\ref{10}.


\begin{thebibliography}{99}
\bibitem{oyz} C. Ouyang, W. Yang, and R. Zhao, Characterizations of
Bergman spaces and Bloch space in the unit ball of $\cn$, {\it Trans.
Amer. Math. Soc.} {\bf 347} (1995), 4301-4313.
\bibitem{rudin} W. Rudin, {\it Function Theory in the Unit Ball of $\cn$},
Springer-Verlag, New York, 1980.
\bibitem{zhu} K. Zhu, {\it Spaces of Holomorphic Functions in the Unit Ball},
Springer-Verlag, New York, 2004.
\end{thebibliography}
\end{document}